\documentclass[12pt]{article}
\usepackage[amsmath]{e-jc}
\usepackage{enumitem}

\usepackage{graphicx}
\usepackage{subcaption}
\usepackage{listings}
\usepackage{lipsum}

\theoremstyle{plain}
\newtheorem{result}[theorem]{Result}
\numberwithin{theorem}{section}

\dateline{April 12, 2023}{TBD}{TBD}


\Copyright{The authors. Released under the CC BY-ND license (International 4.0).}

\title{Improved lower bounds for Queen's Domination \\ via an exactly-solvable relaxation}

\author{Archit Karandikar\footnote{Works at Waymo LLC.}\\
\small\tt architkarandikar@gmail.com
\and
Akashnil Dutta\footnote{Works at LinkedIn Inc.}\\
\small\tt akashnil07@gmail.com}

\begin{document}

\maketitle

\begin{abstract}
  The Queen's Domination problem, studied for over 160 years, poses the following question: What is the least number of queens that can be arranged on a $m \times n$ chessboard so that they either attack or occupy every cell? 
  
  We propose a novel relaxation of the Queen's Domination problem and show that it is exactly solvable on both square and rectangular chessboards. As a consequence, we improve on the best known lower bound for rectangular chessboards in $\approx 12.5\%$ of the non-trivial cases. As another consequence, we simplify and generalize the proofs for the best known lower-bounds for Queen's Domination of square $n \times n$ chessboards for $n \equiv \{0,1,2\} \mod 4$ using an elegant idea based on a convex hull.
  
  Finally, we show some results and make some conjectures towards the goal of simplifying the long complicated proof for the best known lower-bound for square boards when $n \equiv 3 \mod 4$ (and $n > 11$). These simple-to-state conjectures may also be of independent interest.
\end{abstract}


\section{Introduction}

The Queen's graph $Q_{m \times n}$ has the squares of an $m \times n$ chessboard as its vertex set $V(Q_{m \times n})$, with two vertices connected by an edge if they are in the same row, column or diagonal. A subset $\mathcal D$ of the vertex set is said to dominate $Q_{m \times n}$ if every vertex is either contained in or is adjacent to an element in $\mathcal D$. The problem of finding the smallest dominating set of $Q_{m \times n}$  - the size of which is denoted by $\gamma(Q_{m \times n})$ - is the Queen's Domination problem.

We identify the $mn$ squares of the $m \times n$ chessboard by their co-ordinates $(x,y)$ where $x \in \{0,...,m-1\}$ denotes column number (increasing left to right) and $y \in \{0,...,n-1\}$ denotes row number (increasing bottom to top). Diagonals which go left as they go up are called sum diagonals since the value $x+y$ is fixed long them. Diagonals which go right as they go up are called difference diagonals and the value $x-y$ is fixed along them.

The Queen's Domination problem is a special case of the set cover problem, which is itself a discrete constrained optimization problem. It can be stated as follows:

\begin{definition}[Queen's Domination Problem]\label{Def:QueensDomination}
\begin{align*}
\text{$\gamma(Q_{m \times n})=$ minimize } & p \\ 
\text{such that } & \text{$\exists$ rows } R_i\text{, }0\leq i\leq p-1 \\
& \text{$\exists$ columns } C_i\text{, }0\leq i\leq p-1 \\
& \text{$\exists$ sum diagonals } S_i\text{, }0\leq i\leq p-1 \\
& \text{$\exists$ difference diagonals } D_i\text{, }0\leq i\leq p-1 \\
\text{subject to } & V(Q_{m \times n}) = (\bigcup_{i=0}^{p-1}R_i) \cup (\bigcup_{i=0}^{p-1}C_i) \cup (\bigcup_{i=0}^{p-1}S_i) \cup (\bigcup_{i=0}^{p-1}D_i) \\
\text{and }& R_i \cap C_i \cap S_i \cap D_i \neq \phi\text{, }0\leq i\leq p-1
\end{align*}
\end{definition}

The row $R_i$, column $C_i$ and diagonals $S_i, D_i$ are the ones occupied by the $i^\text{th}$ queen. The last constraint in the above formulation, which requires that the $i^\text{th}$ row, column, sum diagonal and difference diagonal all intersect at a square on the chessboard, corresponds to queen placement and seems a little hard to tackle. With this motivation, we relax this constraint and refer to the resultant problem as the Relaxed Queen's Domination problem. The relaxation - whose solution we denote by $\beta(Q_{m\times n})$ - asks for the least number $p$ such that $p$ rows, $p$ columns, $p$ sum diagonals and $p$ difference diagonals suffice to cover a $m \times n$ chessboard.

\begin{definition}[Relaxed Queen's Domination Problem]\label{Def:RelaxedQueensDomination}
\begin{align*}
\text{$\beta(Q_{m \times n})=$ minimize } & p \\
\text{such that } & \text{$\exists$ rows } R_i\text{, }0\leq i\leq p-1 \\
& \text{$\exists$ columns } C_i\text{, }0\leq i\leq p-1 \\
& \text{$\exists$ sum diagonals } S_i\text{, }0\leq i\leq p-1 \\
& \text{$\exists$ difference diagonals } D_i\text{, }0\leq i\leq p-1 \\
\text{subject to } & V(Q_{m \times n}) = (\bigcup_{i=0}^{p-1}R_i) \cup (\bigcup_{i=0}^{p-1}C_i) \cup (\bigcup_{i=0}^{p-1}S_i) \cup (\bigcup_{i=0}^{p-1}D_i) \stepcounter{equation}\tag{\theequation}\label{eqn:CoverTheBoard}
\end{align*}
\end{definition}

Note that since we've relaxed a constraint, we know that relaxed Queen's Domination is a lower bound on Queen's Domination i.e. $\beta(Q_{m \times n}) \leq \gamma(Q_{m \times n})$. This relaxation is the crux of our work. It is remarkable that this relaxation is both (1) easy enough to be exactly solvable using simple constructions and impossibility proofs (2) hard enough to improve the best-known lower bound in a few cases and match it in most cases. We will solve this relaxation exactly for both square and rectangular boards, simplifying the proofs for the former's lower bounds and improving the latter's lower bounds in the process. As a consequence, \textit{throughout this work, we will deal with choosing queen-occupied rows, columns and diagonals rather than placing queens.}

\section{Related Work}

The Queen's Domination problem on square boards was included (Problem C18) in a collection of unsolved problems \cite{UnsolvedProblemNTBook} in Number Theory in 1994. It was first considered by de Jaenisch \cite{jaenisch1862} in 1862, when he gave minimum dominating sets of $Q_{n \times n}$ for $n\leq 8$. Ball \cite{ball1892mathematical} considered several new questions about Queen's Domination. Ahrens \cite{Ahrens1901} and von Szily \cite{Szily1902, Szily1903} gave minimum dominating sets of $Q_{n \times n}$ for $9\leq n\leq 13$ and $n=17$. Many of these were shown to be minimum by recent works on lower bounds. These works also made efforts to list all minimum dominating sets for each $n$, producing lists modulo symmetries. Ahrens \cite{Ahrens1910} summarized these results in his 1910 book.

The first non-trivial lower bound for square boards was obtained by Raghavan and Venkatesan in 1987 \cite{Raghavan1987281} when they showed that $\gamma(Q_{n \times n}) \geq (n-1)/2$. Another proof of this bound was given by Spencer in 1990 (as cited and stated in \cite{Cockayne199013, Weakley1995, BurgerThesis}). Weakley \cite{Weakley1995} improved the bound for case $n \equiv 1 \mod 4$, showing that $\gamma(Q_{n \times n}) \geq (n+1)/2$ in this case. Watkins \cite{Watkins2004} collected these results in his comprehensive treatise of mathematical chessboard problems. {\"O}sterg{\aa}rd and Weakley \cite{stergrd2001ValuesOD} showed that these established bounds were all very good (either exact or off-by-1) for $n \leq 120$ by adapting an algorithm from Knuth \cite{Knuth2000} originally developed for the exact cover problem. Weakley and Finozhenok \cite{Weakley2002, Finozhenok2007, Weakley1995} improved the bound for the case  $n \equiv 3 \mod 4$ (when $n > 11$) to $(n+1)/2$ with a long complex proof spanning 3 papers and more than 35 pages. For even square boards, there has been no direct improvement over the very first $(n-1)/2$ lower bound, but Weakley \cite{Weakley2022} showed several results and made some conjectures towards this goal, especially for the case $n \equiv 0 \mod 4$.

The progression of the best-known lower bound for square boards over the years is shown in Table \ref{Table:SquareBoardsSummary}. In this work, we simplify the proofs for all these lower bounds, except for this last bound for the case  $n \equiv 3 \mod 4$ (when $n > 11$) which has the aforementioned long complex proof. In particular, we simplify the proof for the best-known lower bound for $n \equiv \{0,1,2\} \mod 4$. We make some conjectures towards the goal of simplifying this last bound which may be of independent interest.

\begin{table}[h!]
\centering
\begin{tabular}{|l||l|l|l||l|} 
 \hline
 Case & \textbf{\shortstack{Raghavan et al \\ \cite{Raghavan1987281}, Spencer \cite{Cockayne199013}}} & \textbf{Weakley \cite{Weakley1995}} & \shortstack{Weakley et al \\ \cite{Weakley2002, Finozhenok2007, Weakley1995}} & Our Results\\
 \hline \hline
 \textbf{n=4k} & \textbf{$\gamma(Q_n)\geq 2k$} & \textbf{$\gamma(Q_n)\geq 2k$} & \textbf{$\gamma(Q_n)\geq 2k$} & $\beta(Q_n)=2k$ \\
 \hline
 \textbf{n=4k+1} & \textbf{$\gamma(Q_n)\geq 2k$} & \textbf{$\gamma(Q_n)\geq 2k+1$} & \textbf{$\gamma(Q_n)\geq 2k+1$} & $\beta(Q_n)=2k+1$ \\ 
 \hline
 \textbf{n=4k+2} & \textbf{$\gamma(Q_n)\geq 2k+1$} & \textbf{$\gamma(Q_n)\geq 2k+1$} & \textbf{$\gamma(Q_n)\geq 2k+1$} & $\beta(Q_n)=2k+1$ \\
 \hline
 n=4k+3 & \textbf{$\gamma(Q_n)\geq 2k+1$} & \textbf{$\gamma(Q_n)\geq 2k+1$} & \shortstack{$\gamma(Q_n)\geq 2k+2$ \\ for $n>11$} & $\beta(Q_n)=2k+1$ \\ 
 \hline
 \hline
 Year & 1987 & 1995 & 2007 & 2023 \\
 \hline
\end{tabular}
\caption{Progression of the lower bound for square boards over the years. Note that $Q_n$ is used in this table as a shorthand for $Q_{n \times n}$. The cells highlighted in bold are those for whose results we simplify the proof in this work. Since $\beta(Q_n) \leq \gamma(Q_n)$, our work not only simplifies but also generalizes the results in these highlighted cells.}
\label{Table:SquareBoardsSummary}
\end{table}

The first lower bound for the general case of rectangular boards was obtained by Raghavan and Venkatesan in 1987 \cite{Raghavan1987281} when they showed that $\gamma(Q_{m \times n})\geq \min\{m-1,n-1,(m+n-2)/4\}$. We refer to the boards where $\max\{m,n\} \geq 3\min\{m,n\}-2$ as trivial boards. Bozoki et al \cite{Bozoki2019} marginally improved the lower bound to $\min\{m,n\}$ on the trivial case and $(m+n-2)/4$ on the non-trivial case. They posed the question of whether this bound can be improved which we answer in the affirmative. In their work they also computed the exact values of $\gamma(Q_{m \times n})$ for $1\leq m,n\leq 18$ and showed some other results.

\begin{figure}[t]
\includegraphics[width=16cm]{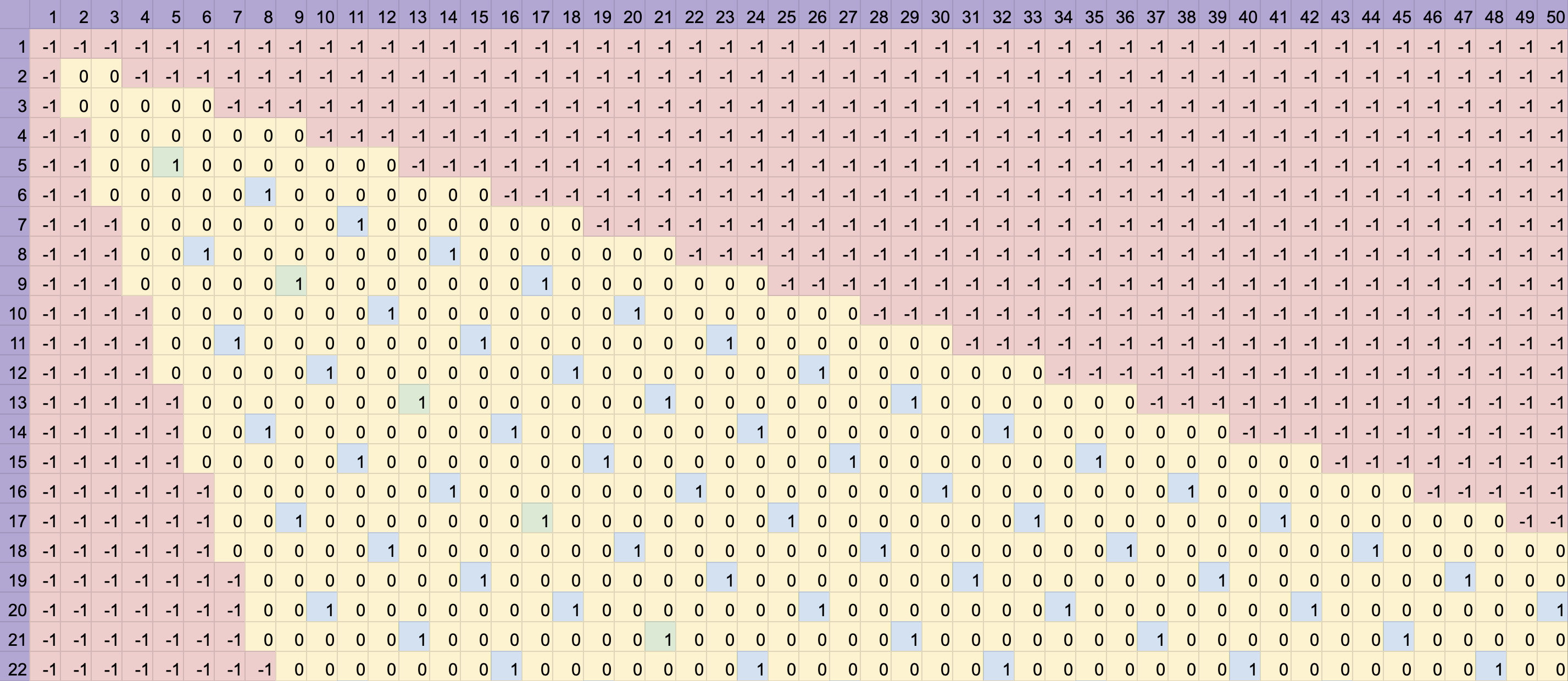}
\centering
\caption{Summary of our improvements to the lower bound for $\gamma(Q_{m \times n})$. Row and column indices denote values of $m$ and $n$ respectively. A value of $-1$ (red) indicates boards on which $\gamma(Q_{m \times n})$ is trivially equal to $\min\{m,n\}$. A value of $1$ (blue) indicates the boards on which we improve upon the $(m+n-2)/4$ lower bound by $1$. The improvement for square boards (green) was already established by Weakley \cite{Weakley1995}. The blue cells account for $\sim 12.5\%$ of the non-trivial cases since one in every 8 non-trivial cases in every row/column isn't coloured yellow. A value of $0$ (yellow) indicates the boards on which we match the existing $(m+n-2)/4$ lower bound. Note that on some of these specific small values of $(m,n)$ a better computational lower bound than the generally applicable $(m+n-2)/4$ lower bound is known. This figure does not reflect these computationally known bounds.}
\label{figure:RectangularBoardsSummary}
\end{figure}

In this work, we improve the best known lower bound for $\approx 12.5\%$ of the non-trivial boards, as shown in Figure \ref{figure:RectangularBoardsSummary}. The essence of this work can be found in Definition \ref{Def:RelaxedQueensDomination}, Lemma \ref{Lemma:SpacedOddRectangularGrid} and Theorem \ref{Thm:4kPlus1RectangularBoardsGeneralization}. In Theorem \ref{Thm:RectangularBoardsFullSolution} we show that \[
\gamma(Q_{m \times n})\geq\beta(Q_{m \times n})= \begin{cases}
\min\{m,n\} & \text{if $\max\{m,n\} \geq 3\min\{m,n\}-2$}\\
(m+n-2)/4 + 1 & \text{else if $m,n$ are even and $m+n \equiv 6 \mod 8$}\\
(m+n-2)/4 + 1 & \text{else if $m,n$ are odd and $m+n \equiv 2\mod 8$}\\
\lceil (m+n-2)/4 \rceil & \text{otherwise}
 \end{cases}
\]

\section{Relaxed Queen's Domination on Square Boards}

In this section we will completely solve the Relaxed Queen's Domination problem on square boards (Theorem \ref{Thm:SquareBoardsFullSolution}). We've already noted that this will yield a lower bound on the Queen's Domination problem since $\gamma(Q_{n\times n}) \geq \beta(Q_{n\times n})$. As a result, we will both generalize the best known lower bounds on $\gamma(Q_{n\times n})$ for $n \equiv \{0,1,2\} \mod 4$ and simplify their proofs in Thoerems \ref{Thm:AllSquareBoards} and \ref{Thm:4k1SquareBoards}. In doing so, we will also completely solve the Relaxed Bishop's Domination problem (Definition \ref{Def:RelaxedBishopsDomination}) on square boards in Lemma \ref{Lemma:BishopsSquareBoardsFullSolution}.

\begin{theorem}\label{Thm:AllSquareBoards}
    For all natural numbers $n$, \[\gamma(Q_{n\times n}) \geq \beta(Q_{n\times n}) \geq (n-1)/2\]
\end{theorem}
\begin{proof}
  Let if possible for some $p<(n-1)/2$, there be a choice of $p$ rows, $p$ columns, $p$ sum diagonals and $p$ difference diagonals which satisfy Equation \ref{eqn:CoverTheBoard} to cover $V(Q_{n \times n})$. Since $p\geq 1$, we must have $n \geq 4$ and $n-p > 2$.
  
  Amongst the set $U$ of $(n-p)\times (n-p)$ cells uncovered by the chosen rows and columns, consider the subset $K$ of the outermost ones - those which touch the boundary of the convex hull of $U$. There are $4(n-p)-4$ points in $K$ (since $n-p\geq 2$) and each of the $2p$ diagonals can cover at most $2$ of these. This implies that \begin{align*}
    2 \times 2p &\geq 4(n-p)-4 \stepcounter{equation}\tag{\theequation}\label{eqn:AllSquareBoardsEquation} \\
    p &\geq (n-1)/2 
  \end{align*}
\end{proof}

\begin{figure}[t]
\includegraphics[width=16cm]{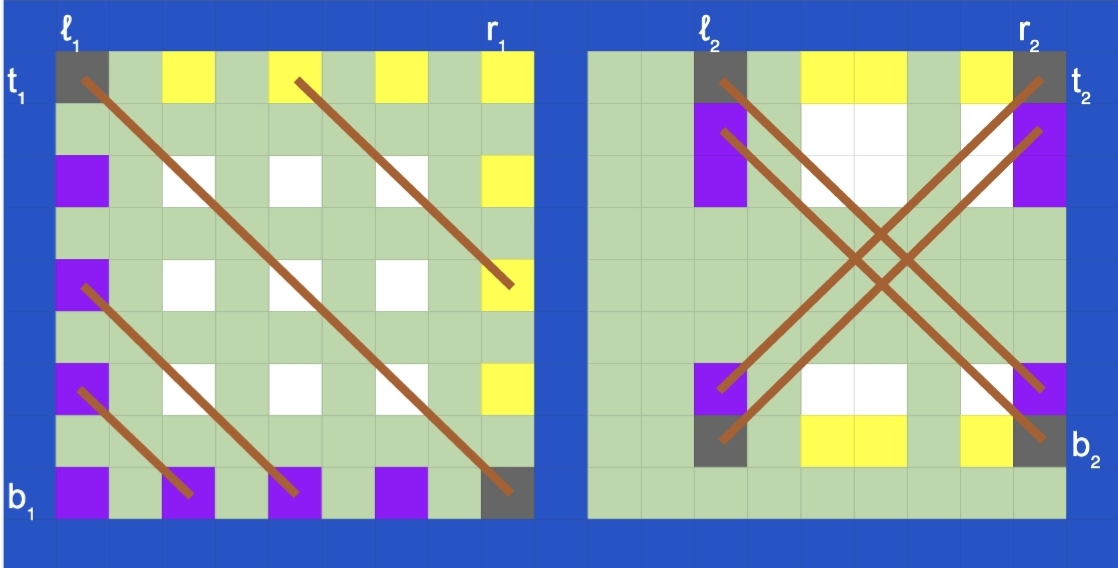}
\centering
\caption{Illustration of the proofs of Theorems \ref{Thm:AllSquareBoards} and \ref{Thm:4k1SquareBoards}. In both boards, we have $n=4k+1=9$, $k=2$ and $p=2k=4$. The blue boundary isn't part of either board. The chosen rows $(R_i)_{i=1}^p$ and columns $(C_i)_{i=1}^p$ are coloured green.
Cells not coloured green constitute the set $U$. Cells coloured grey, yellow($K_Y$) and purple($K_P$) constitute the set $K$ of the outermost cells of $U$. 
The boards on the left and right depict cases 1 and 2 in the proof of Theorem \ref{Thm:4k1SquareBoards} respectively.}
\label{figure:SquareBoardsProof}
\end{figure}

This proof of Theorem \ref{Thm:AllSquareBoards} (illustrated in Figure \ref{figure:SquareBoardsProof}) both generalizes and simplifies the proofs that were obtained by Raghavan and Venkatesan in 1987 \cite{Raghavan1987281} and Spencer in 1990 \cite{Cockayne199013, Weakley1995, BurgerThesis}. While these previously known proofs weren't complex, this one is the simplest because the set of cells $K$ defined above (based on the convex hull of $U$) precisely captures the difficulty of choosing diagonals to cover the cells $U$ uncovered by the chosen rows and columns. The only work we know which specifically considers this set $K$ is \cite{Weakley2022}. Theorem \ref{Thm:AllRectangularBoards} in the next section generalizes this result further for rectangular boards.

\begin{theorem}\label{Thm:4k1SquareBoards}
  For all natural numbers $n$ of the form $4k+1$ where $k$ is an integer, \[\gamma(Q_{n\times n}) \geq \beta(Q_{n\times n}) \geq (n+1)/2\]
\end{theorem}
\begin{proof}
Consider $n$ of the form $4k+1$. Building on the result from Theorem \ref{Thm:AllSquareBoards}, all we need to show is that $\beta(Q_{n\times n})$ cannot equal $(n-1)/2=2k$. Let if possible $\beta(Q_{n\times n})=2k=p$, so that there are $p$ rows, $p$ columns, $p$ sum diagonals and $p$ difference diagonals which satisfy Equation \ref{eqn:CoverTheBoard} to cover $V(Q_{n \times n})$. This can only happen if equality holds in Equation \ref{eqn:AllSquareBoardsEquation}; each of the $2p$ diagonals covers exactly two cells in $K$ and these $2p$ pairs partition $K$.

Consider again the subset $K$ of outermost cells amongst the set $U$ of $(n-p)\times (n-p)$ cells uncovered by the chosen rows and columns. Let $l$, $r$ denote the leftmost, rightmost unchosen rows and $b$, $t$ denote the bottommost, topmost unchosen columns.

\textbf{\textit{Case 1, $r-l=t-b$:}} The cells $(r,b)$ and $(l,t)$ must be connected by a sum diagonal $S^*$ since otherwise the two difference diagonals covering them each cover only one cell in $K$, an impossibility. We colour the $K$-cells $K_Y=\{(x,y)\in K\ | x+y>r+b\}$ above $S^*$ yellow and those $K_P=\{(x,y)\in K\ | x+y<r+b\}$ below $S^*$ purple.

Note that $|K_Y|=|K_P|=4k-1$. Each of the $p=2k$ chosen difference diagonals covers one yellow and one purple cell. Each of the $2k-1$ chosen sum diagonals besides $S^*$ covers either two yellow (if above $S^*$) or two purple cells (if below $S^*$). Since every $K$-cell is covered exactly once, the number of sum diagonals above $S^*$ and below $S^*$ must be equal, which is impossible since $2k-1$ is odd.

\textbf{\textit{Case 2, $r-l\neq t-b$:}} WLG, $r-l < t-b$. We colour the $K_Y=\{(x,y)\in K | y\in\{b,t\}, l<x<r\}$ top and bottom $K$-rows  (besides corners) yellow and the $K_P=\{(x,y)\in K | x\in\{l,r\}, b<y<t\}$ left and right $K$-columns (besides corners) purple.

Note that $|K_Y|=|K_P|=4k-2$. The four diagonals which cover the four corners must each also cover an additional purple cell since $r-l < t-b$. There are four more yellow cells ($4k-2$) than purple cells ($4k-6$) left to cover using the remaining $4k-4$ diagonals. Thus, at least one of the diagonals must cover two yellow cells, which contradicts $r-l < t-b$.

\end{proof}

As before, this proof of Theorem \ref{Thm:4k1SquareBoards} (illustrated in Figure \ref{figure:SquareBoardsProof}) both generalizes and simplifies the proof by Weakley from 1995 \cite{Weakley1995}. The reason this proof is simpler is the same as before: the subset $K$ of outermost cells of the set of cells $U$ uncovered by chosen rows and columns precisely captures the difficulty of covering them using chosen diagonals.

We will now completely solve the Relaxed Queen's domination problem in Theorem \ref{Thm:SquareBoardsFullSolution} by showing that these lower bounds are also upper bounds. Towards that end, we first define and completely solve the Relaxed Bishop's Domination problem. $B_{m \times n}$ is a graph which has the squares of an $m \times n$ chessboard as its vertex set $V(B_{m \times n})$, with two vertices connected by an edge if they are in the same diagonal. The problem of finding the smallest dominating set of $B_{m \times n}$ is the Bishop's Domination problem. This is a solved problem \cite{Cockayne1986, Watkins2004} and we are interested only in it's relaxation which is defined as follows:

\begin{definition}[Relaxed Bishop's Domination Problem]\label{Def:RelaxedBishopsDomination}
\begin{align*}
\text{$\alpha(B_{m \times n})=$ minimize } & p \\
\text{such that } & \text{$\exists$ sum diagonals } S_i, 1\leq i\leq p \\
& \text{$\exists$ difference diagonals } D_i, 1\leq i\leq p \\
\text{subject to } & V(B_{m \times n}) = (\bigcup_{i=1}^pS_i) \cup (\bigcup_{i=1}^pD_i) \stepcounter{equation}\tag{\theequation}\label{eqn:BishopsCoverTheBoard}
\end{align*}
\end{definition}

A useful simple observation is that the Relaxed Queen's and Bishop's Domination problems are trivially monotonic. Monotonicity doesn't hold for the Queen's Domination problem since $\lambda(Q_{8\times 11})=6>5=\lambda(Q_{9\times 11})=\lambda(Q_{10\times 11})=\lambda(Q_{11\times 11})$ as shown in \cite{Weakley2018}. 

\begin{observation}\label{Obs:Monotonicity}
If $m\leq m'$, $n\leq n'$ then $\alpha(B_{m \times n}) \leq \alpha(B_{m' \times n'})$, $\beta(B_{m \times n}) \leq \beta(B_{m' \times n'})$.
\end{observation}

\begin{lemma}\label{Lemma:BishopsSquareBoardsFullSolution}
  For all natural numbers $n$
  \[
    \alpha(B_{n \times n})= \begin{cases}
    2k-1 & \text{if $n=2k$}\\
    2k+1 & \text{if $n=2k+1$}
     \end{cases}
  \]
\end{lemma}
\begin{proof}
The statement is equivalent to $\alpha(B_{n \times n})=2k-1$ if $n\in\{2k-1,2k\}$. Given \ref{Obs:Monotonicity}, it suffices to show that $\alpha(B_{(2k-1) \times (2k-1)}) \geq 2k-1$ and $\alpha(B_{{2k} \times {2k}}) \leq 2k-1$.

Let if possible $\alpha(B_{(2k-1) \times (2k-1)})<2k-1$. Append $2k-2$ rows to the right and $2k-2$ columns to the top of $V(B_{(2k-1) \times (2k-1)})$ and choose these in the solution to the Relaxed Queen's Domination for this larger $V(Q_{(4k-3)\times(4k-3)})$ board. The cells uncovered by the chosen rows and columns form the bottom-left $V(B_{(2k-1) \times (2k-1)})$ sub-board that we started with, which can be covered by at most $2k-2$ sum and difference diagonals since $\alpha(B_{(2k-1) \times (2k-1)})<2k-1$. It follows that $\beta(Q_{(4k-3) \times (4k-3)}) \leq 2k-2$, which contradicts Theorem \ref{Thm:4k1SquareBoards}. Hence, we have shown that $\alpha(B_{(2k-1) \times (2k-1)})\geq 2k-1$.

Now consider the $V(B_{2k\times 2k})$ board. The set of possible sums $\mathcal S = \{x+y|(x,y)\in V(B_{2k\times 2k})\}=\{0,\ldots,4k-2\}$ and the set of possible differences $\mathcal D = \{x+y|(x,y)\in V(B_{2k\times 2k})\}=\{-(2k-1),\ldots,2k-1\}$. We choose the set of sum diagonals with odd sums $\{s \in \mathcal S|s\text{ is odd}\}$ and the set of difference diagonals with even differences $\{d \in \mathcal S|d\text{ is even}\}$, both of size $2k-1$. Together these cover $V(B_{2k\times 2k})$ since for each $(x,y)$ the parity of $x+y$ is the same as the parity of $x-y$. It follows that $\alpha(B_{2k\times 2k})\leq 2k-1$.
\end{proof}

\begin{theorem}\label{Thm:SquareBoardsFullSolution}
  For all natural numbers $n$
  \[
    \beta(Q_{n \times n})= \begin{cases}
    2k & \text{if $n=4k$}\\
    2k+1 & \text{if $n=4k+1$}\\
    2k+1 & \text{if $n=4k+2$}\\
    2k+1 & \text{if $n=4k+3$}
     \end{cases}
  \]
\end{theorem}
\begin{proof} Consider $n\in\{4k,4k+1,4k+2,4k+3\}$ and define $f(n)=2k$ for $n=4k$ and $f(n)=2k+1$ otherwise. Theorems \ref{Thm:AllSquareBoards} and \ref{Thm:4k1SquareBoards} have already shown that $\beta(Q_{n \times n}) \geq f(n)$.

We choose the rightmost $f(n)$ rows and the topmost $f(n)$ columns of $V(Q_{n \times n})$. Let $U$ denote the cells uncovered by the chosen rows and columns. If $n=4k$, $U=V(B_{2k \times 2k})$ which we know from Lemma \ref{Lemma:BishopsSquareBoardsFullSolution} can be covered with $2k-1<f(n)$ sum and difference diagonals, implying that $\beta(Q_{n \times n}) \leq f(n)$. If $n=4k+3$, $U=V(B_{(2k+2) \times (2k+2)})$ which we know from Lemma \ref{Lemma:BishopsSquareBoardsFullSolution} can be covered with $2k+1=f(n)$ sum and difference diagonals. Observation \ref{Obs:Monotonicity} implies that $\beta(Q_{n \times n}) \leq f(n)=2k+1$ for $n \in \{4k+1,4k+2,4k+3\}$.
\end{proof}

\begin{figure}[t]
\includegraphics[width=16cm]{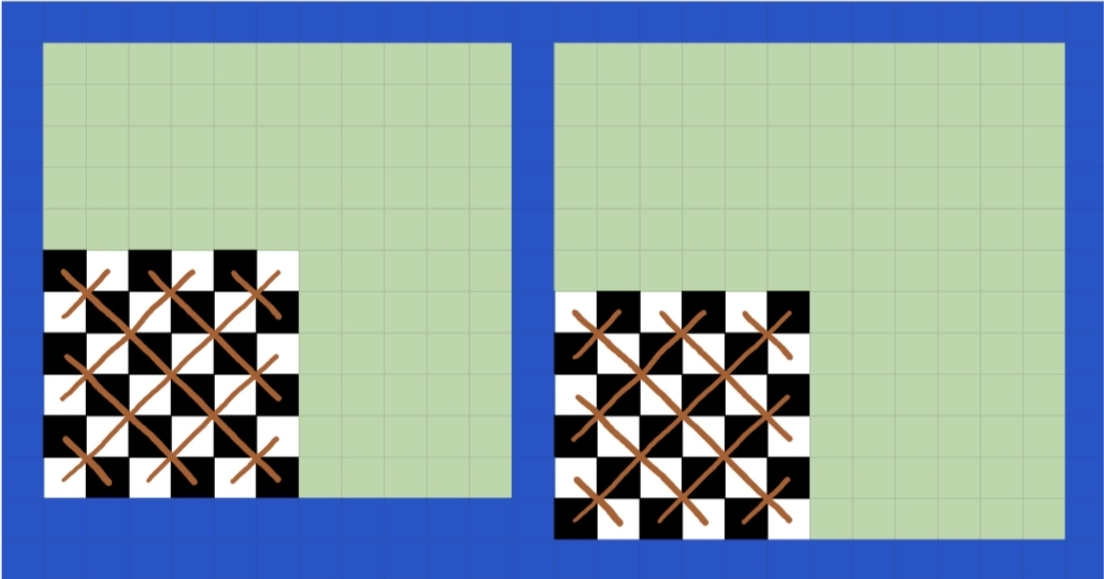}
\centering
\caption{Illustration of the proof of Theorem \ref{Thm:SquareBoardsFullSolution}. Blue cells aren't part of either board. In the first board on the left we have $n_1=4k_1+3=11,f(n_1)=5$ and in the second we have $n_2=4k_2=12,f(n_2)=6$. The chosen rows $(R_i)_{i=1}^p$ and columns $(C_i)_{i=1}^p$ are coloured green and the chosen diagonals are shown in brown. The one extra remaining sum and difference diagonal in the figure on the right isn't shown and can be chosen arbitrarily. The black-and-white section of the two boards also illustrates the chosen diagonals in Lemma \ref{Lemma:BishopsSquareBoardsFullSolution} for $n_3=2k_3=6$.}
\label{figure:SquareBoardsConstruction}
\end{figure}

This proof of Theorem \ref{Thm:SquareBoardsFullSolution} (illustrated in Figure \ref{figure:SquareBoardsConstruction}) completes our solution of the Relaxed Queen's Domination problem for square boards.

\section{Relaxed Queen's Domination on Rectangular Boards}

In this section we will completely solve the Relaxed Queen's Domination problem on rectangular boards (Theorem \ref{Thm:RectangularBoardsFullSolution}) and hence improve the lower bound for the Queen's Domination problem on $\approx 12.5\%$ of the non-trivial rectangular boards as shown in Figure \ref{figure:RectangularBoardsSummary}. Moreover, we generalize and simplify the proof of the currently known lower bound for Queen's Domination of rectangular boards (Theorem \ref{Thm:AllRectangularBoards}). Theorem \ref{Thm:4kPlus1RectangularBoardsGeneralization} generalizes Theorem \ref{Thm:4k1SquareBoards} and Weakley's improvement \cite{Weakley1995} for the lower bound on $\gamma(Q_{(4k+1)\times(4k+1)})$ to rectangular boards. It uses Lemma \ref{Lemma:SpacedOddRectangularGrid}, which contains the central idea of this work. We will also completely solve the Relaxed Bishop's Domination problem on rectangular boards in Lemma \ref{Lemma:BishopsRectangularBoardsFullSolution}. This section generalizes the previous section to rectangular boards.

\begin{theorem}\label{Thm:AllRectangularBoards}
    For all natural numbers $m,n$, \[\gamma(Q_{m\times n}) \geq \beta(Q_{m\times n}) \geq 
    \begin{cases}
    \min\{m,n\} & \text{if $\max\{m,n\} \geq 3\min\{m,n\}-2$}\\
    (m+n-2)/4 & \text{otherwise}
     \end{cases}\]
\end{theorem}
\begin{proof}
Let $m \geq n$ WLG. Let $\beta(Q_{m\times n})=p$ so that there is a choice of $p$ rows, $p$ columns, $p$ sum diagonals and $p$ difference diagonals which satisfy Equation \ref{eqn:CoverTheBoard} to cover $V(Q_{m \times n})$.

In the trivial case where $m\geq 3n-2$, let if possible $p<n$. Consider an unchosen row of $m$ cells. Each of the other $3p \leq 3(n-1)$ chosen lines (columns and diagonals) can cover at most one cell in it. Hence $m\leq3(n-1)<3n-2\leq m$, a contradiction.

For the non-trivial case note that $m<3n-2 \Leftrightarrow (m+n-2)/4 < n-1$. Let if possible $p<(m+n-2)/4<n-1$. Amongst the set $U$ of $(m-p)\times (n-p)$ cells uncovered by the chosen rows and columns, consider again the subset $K$ of the outermost ones - those which touch the boundary of the convex hull of $U$. There are $2(m-p)+2(n-p)-4$ points in $K$ (since $m-p \geq n-p> 1$) and each of the $2p$ diagonals can cover at most $2$ of these. This implies that \begin{align*}
    2 \times 2p &\geq 2(m-p)+2(n-p)-4 \\
    p &\geq (m+n-2)/4 
\end{align*} 
\end{proof}

\begin{figure}[t]
\includegraphics[width=12cm]{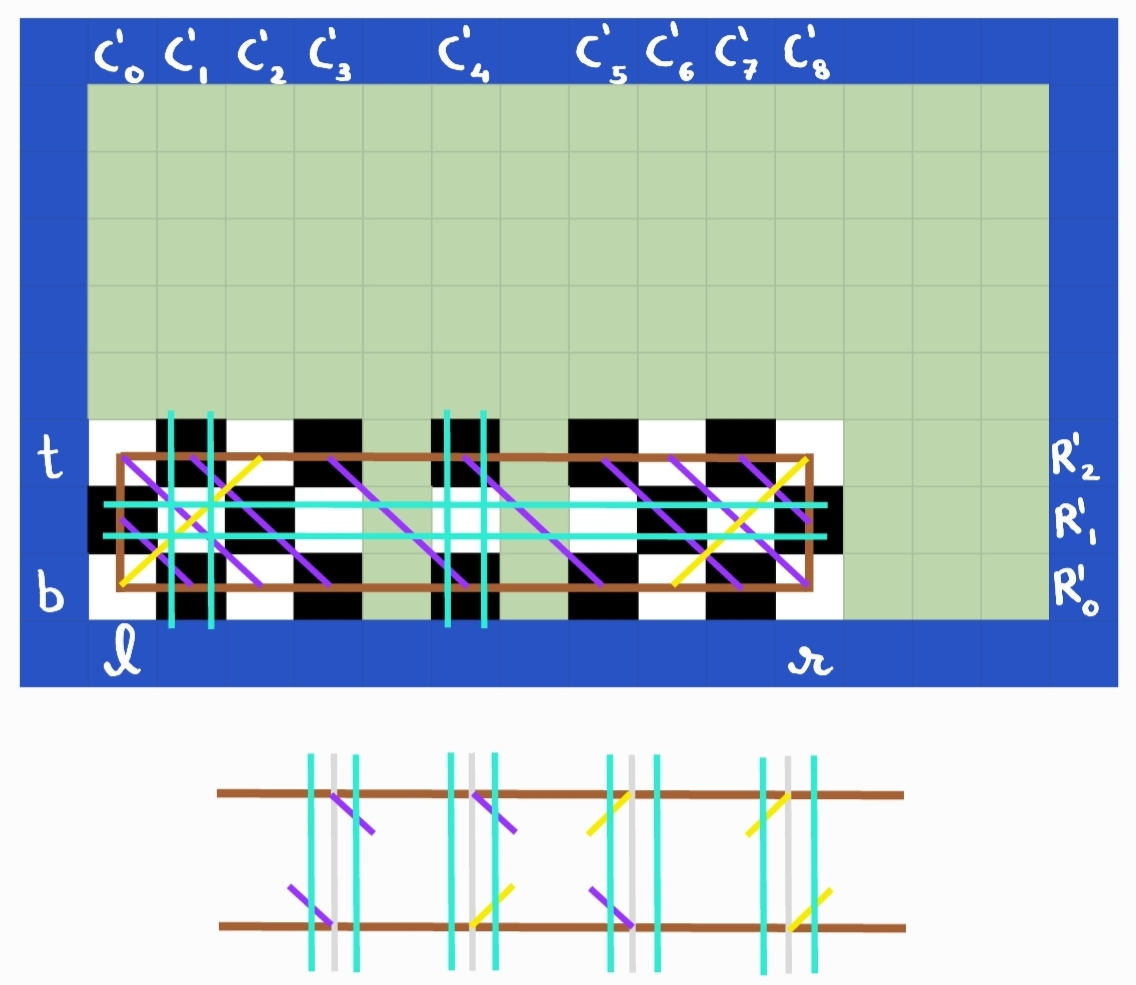}
\centering
\caption{Illustration of the proofs of Lemma \ref{Lemma:SpacedOddRectangularGrid} ($m'=9,n'=3$) and Theorem \ref{Thm:4kPlus1RectangularBoardsGeneralization} ($m=14,n=8$). The black and white cells form the spaced grid $G$. The green cells are those covered by the chosen rows and columns. Rectangle $\mathcal T$ is shown in brown. Purple and yellow segments correspond to chosen sum and difference diagonals respectively. Vertical ($v$) and horizontal ($h$) cross-sections are cyan coloured. We see that $q(v)=1$ and $q(h)=5$ are invariant. As per equation \ref{eqn:PurpleYellowDiagonalDifferenceEquation}, $q(v)(r-l)=q(h)(t-b)=10$  and difference between the total length of purple and yellow segments $l_P-l_Y=10\sqrt 2$. The $4$ cases for the transition of the vertical cross-section are shown at the bottom.}
\label{figure:RectangularBoardsProof}
\end{figure}

This proof of Theorem \ref{Thm:AllRectangularBoards} (illustrated in Figure \ref{figure:RectangularBoardsProof}) both generalizes and simplifies the proofs by Raghavan and Venkatesan \cite{Raghavan1987281} and Bozoki et al \cite{Bozoki2019}. It also generalizes Theorem \ref{Thm:AllSquareBoards} to rectangular boards. 
Theorem \ref{Thm:4kPlus1RectangularBoardsGeneralization} which follows generalizes Theorem \ref{Thm:4k1SquareBoards} and improves the best known lower bound for the Queen's Domination problem on rectangular boards. Towards this end, we now formally define the `set of uncovered cells' ($U$) structure that we've been using. In Definition \ref{Def:SpacedGrid}, it is useful to think of $\{C'_i\}_{i=0}^{m'-1}$ and $\{R'_i\}_{i=0}^{n'-1}$ as the $m'=m-p$ unchosen columns and $n'=n-p$ unchosen rows.

\begin{definition}[Spaced Grid]\label{Def:SpacedGrid}
Given a set of $m'$ columns $\{C'_i
\}_{i=0}^{m'-1}$ and a set of $n'$ rows $\{R'_i
\}_{i=0}^{n'-1}$, we define the spaced grid with these rows and columns as \[G(\{C'_i\}_{i=0}^{m'-1},\{R'_i\}_{i=0}^{n'-1})=(\bigcup_{i=0}^{m'-1}C'_i) \cap (\bigcup_{i=0}^{n'-1}R'_i)\].
\end{definition}

\begin{lemma}\label{Lemma:SpacedOddRectangularGrid} Given natural numbers $m',n'$, a set of $m'$ columns $\{C'_i
\}_{i=0}^{m'-1}$ and a set of $n'$ rows $\{R'_i
\}_{i=0}^{n'-1}$, let $p$ be the least number such that $p$ sum diagonals $\{S'_i\}_{i=0}^{p-1}$ and $p$ difference diagonals $\{D'_i\}_{i=0}^{p-1}$ suffice to cover the spaced grid $G(\{C'_i\}_{i=0}^{m'-1},\{R'_i\}_{i=0}^{n'-1})$. \[p \geq
\begin{cases}
(m'+n'-2)/2+1 & \text{if $m',n'$ are odd} \\
(m'+n'-2)/2 & \text{otherwise}
\end{cases}\]
\begin{proof} The result is trivial if $\min\{m',n'\}=1$. Otherwise, we again consider the set of cells $K$ which touch the boundary of the convex hull of the spaced grid. Note again that $|K|=(2m'+2n'-4)$ and each of the chosen $2p$ diagonals can cover at most $2$ cells in $C$.
\begin{align*}
    2 \times 2p &\geq 2m'+2n'-4 \stepcounter{equation}\tag{\theequation}\label{eqn:SpacedOddRectangularGridEquation} \\
    p &\geq (m'+n'-2)/2 
\end{align*}
Since $(m'+n'-2)/2$ is integral when $m',n'$ are odd, it suffices to show that equality cannot hold in this case. Let if possible $p=(m'+n'-2)/2$, so that each of the $2p$ chosen diagonals covers exactly $2$ cells in $K$. We assume that columns and rows are indexed left-to-right and bottom-to-top. Let $l=C'_0$, $r=C'_{m'-1}$ denote the leftmost, rightmost columns and $b=R'_0$, $t=R'_{n-1}$ denote the bottommost, topmost rows.

Construct the rectangle $\mathcal T$ with the four vertices at the centers of the cells $(b,l)$, $(b,r)$, $(t,l)$ and $(t,r)$. We construct a purple segment corresponding to each chosen sum diagonal $S'_i$ which joins the centers of cells in it and is confined to $\mathcal T$ (has both endpoints on $\mathcal T$). Similarly, we construct an analogously defined yellow segment corresponding to each chosen difference diagonal $D'_i$. Figure \ref{figure:RectangularBoardsProof} illustrates these constructions.

Consider a vertical cross-section $v$ which intersects $\mathcal T$ and does not pass through the center of any cell. Let $q(v)$ denote the 
difference between the number of purple segments and the number of yellow segments it intersects. Observe that $q(v)$ is invariant as we slide $v$ from left to right. The only places that $v$ stops or starts intersecting some segment is when it crosses over the center of any column. In such transitions, there can be at most two segments that it stops or starts intersecting, as illustrated in the four cases shown at the bottom of Figure \ref{figure:RectangularBoardsProof}. In all cases, $q(v)$ does not change. Similarly, observe that $q(h)$ is invariant as we slide a horizontal cross-section which intersects $\mathcal T$ from bottom to top.

We focus on the difference between the total length $l_P$ of purple segments and the total length $l_Y$ of yellow segments. 
Since $q(v)$ (resp. $q(h)$) is invariant, every $(\delta x)$-width vertical (resp. horizontal) cross-section contributes $\sqrt 2q(v)(\delta x)$ (resp. $\sqrt 2q(h)(\delta x)$) to this quantity and hence \[
l_P-l_Y=\sqrt 2q(v)(r-l)=\sqrt 2q(h)(t-b) \stepcounter{equation}\tag{\theequation}\label{eqn:PurpleYellowDiagonalDifferenceEquation}
\]

Both $q(v)$ and $q(h)$ are non-zero since $m'$, $n'$ (which also equal the number of segments that intersect the horizontal and vertical sides of $\mathcal T$) are odd. Hence, either both $q(v)$ and $q(h)$ are positive or both are negative. But $2q(v)+2q(h)$ should be zero since each of the $p$ sum diagonals (purple segments) and $p$ difference diagonals (yellow segments) intersects $\mathcal T$ exactly twice and since two corners each are covered by purple and yellow segments. This is a contradiction, implying that $p\geq (m'+n'-2)/2 + 1$ if $m',n'$ are odd.
\end{proof}
\end{lemma}

\begin{theorem}\label{Thm:4kPlus1RectangularBoardsGeneralization}
    For all natural numbers $m,n$, \[\gamma(Q_{m\times n}) \geq \beta(Q_{m\times n}) \geq 
    \begin{cases}
    \min\{m,n\} & \text{if $\max\{m,n\} \geq 3\min\{m,n\}-2$}\\
    (m+n-2)/4 + 1 & \text{else if $m,n$ are even and $m+n \equiv 6 \mod 8$}\\
    (m+n-2)/4 + 1 & \text{else if $m,n$ are odd and $m+n \equiv 2\mod 8$}\\
    (m+n-2)/4 & \text{otherwise}
     \end{cases}\]
\end{theorem}
\begin{proof}
We need to improve the lower bound from Theorem \ref{Thm:AllRectangularBoards} by $1$ (a) when $m,n$ are even and $m+n \equiv 6 \mod 8$ and (b) when $m,n$ are odd and $m+n \equiv 2 \mod 8$. Consider such a case and let if possible $\beta(Q_{m\times n})=p=(m+n-2)/4$ so that there is a choice of $p$ rows, $p$ columns, $p$ sum diagonals and $p$ difference diagonals which satisfy Equation \ref{eqn:CoverTheBoard} to cover $V(Q_{m \times n})$. $p$ is odd in case (a) and even in case (b). In both cases, the number of uncovered columns $m'=m-p$ and number of uncovered rows $n'=n-p$ are both odd.

Consider the $m' \times n'$ spaced grid $U$ of uncovered cells. The $p$ chosen sum diagonals and $p$ chosen difference diagonals cover $U$. Since $m',n'$ are odd, Lemma \ref{Lemma:SpacedOddRectangularGrid} implies \begin{align*}
p &\geq (m'+n'-2)/2 + 1 \\
&=(m-p+n-p-2)/2 + 1\\
&= (4p - 2p)/2 + 1 \\
&= p+1
\end{align*}
This contradiction completes the proof.
\end{proof}


\begin{lemma}\label{Lemma:BishopsRectangularBoardsFullSolution}
  For all natural numbers $m,n$
  \[
    \alpha(B_{m \times n})= \begin{cases}
    (m+n-2)/2 + 1 & \text{if $m,n$ are odd}\\
    \lceil (m+n-2)/2 \rceil & \text{otherwise}
     \end{cases}
  \]
\end{lemma}
\begin{proof}
The $V(B_{m \times n})$ board is a special case of a spaced grid with consecutive rows and columns. Lemma \ref{Lemma:SpacedOddRectangularGrid} implies that $\alpha(B_{m \times n}) \geq f(m,n)$, where $f(m,n)=(m+n-2)/2+1$ if $m,n$ are odd and $\lceil(m+n-2)/2\rceil$ otherwise.

We will now show that $\alpha(B_{m \times n}) \leq f(m,n)$. Observation \ref{Obs:Monotonicity} implies that it suffices to show this in the case where $m,n$ are both even, since all other cases can be reduced to this case by increasing $m$ or $n$ or both without increasing $f(m,n)$. The set of possible sums $\mathcal S = \{x+y|(x,y)\in V(B_{m\times n})\}=\{0,\ldots,m+n-2\}$ and the set of possible differences $\mathcal D = \{x-y|(x,y)\in V(B_{m\times n})\}=\{-(n-1),\ldots,m-1\}$. We choose the set of sum diagonals with odd sums $\{s \in \mathcal S|s\text{ is odd}\}$ and the set of difference diagonals with even differences $\{d \in \mathcal D|d\text{ is even}\}$. Each of these is of size $(m+n-2)/2$ since $m,n$ are both even. Together these cover $V(B_{m\times n})$ since for each $(x,y)$, the parity of $x+y$ is the same as the parity of $x-y$. Hence $\alpha(B_{m \times n}) \leq (m+n-2)/2=f(m,n)$.
\end{proof}

\begin{figure}[t]
\includegraphics[width=16cm]{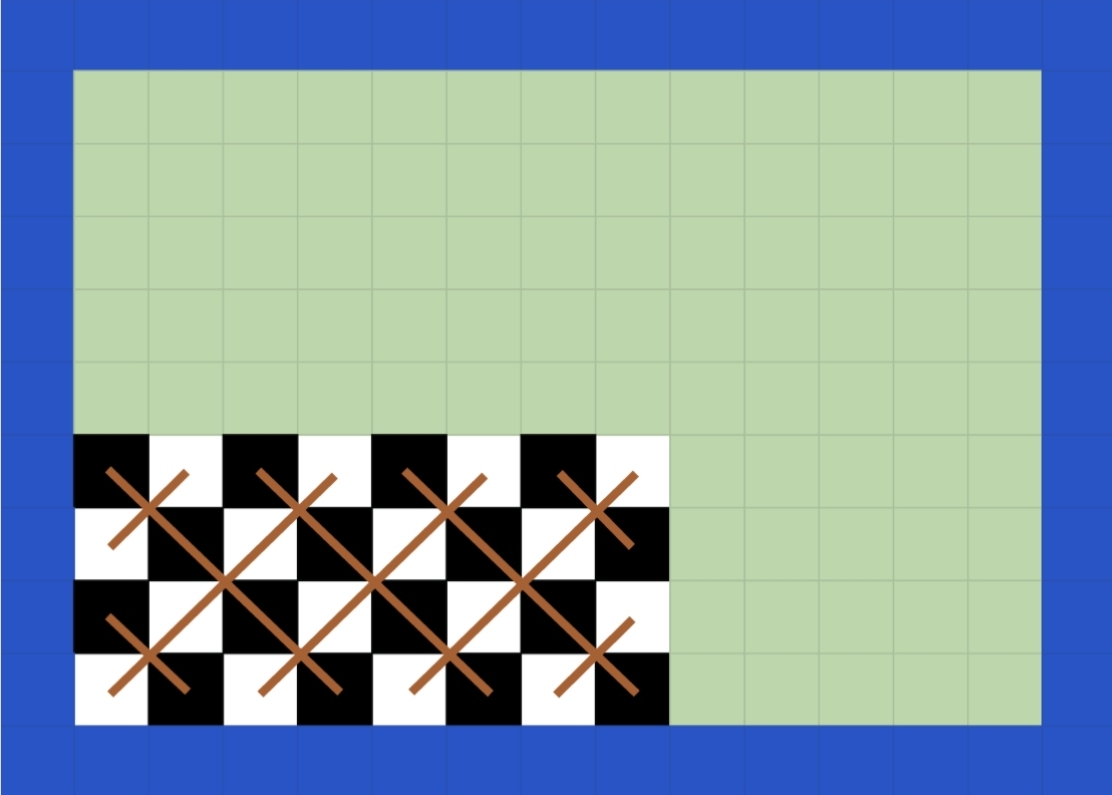}
\centering
\caption{Illustration of the proof of Theorem \ref{Thm:RectangularBoardsFullSolution} for an easy critical board ($m=13$ and $n=9$). Blue cells aren't part of the board. In this case $f(m,n)=(m+n-2)/2=5$. The cells in  $5$ chosen rows and $5$ columns are coloured green. The $5$ chosen sum diagonals $5$ chosen difference diagonals are shown in brown.} 
\label{figure:RectangularBoardsConstruction}
\end{figure}

\begin{theorem}\label{Thm:RectangularBoardsFullSolution}
  For all natural numbers $m,n$
  \[
    \beta(Q_{m \times n})= \begin{cases}
    \min\{m,n\} & \text{if $\max\{m,n\} \geq 3\min\{m,n\}-2$}\\
    (m+n-2)/4 + 1 & \text{else if $m,n$ are even and $m+n \equiv 6 \mod 8$}\\
    (m+n-2)/4 + 1 & \text{else if $m,n$ are odd and $m+n \equiv 2\mod 8$}\\
    \lceil (m+n-2)/4 \rceil & \text{otherwise}
     \end{cases}
  \]
\end{theorem}
\begin{proof}
Theorem \ref{Thm:4kPlus1RectangularBoardsGeneralization} states that $\beta(Q_{m \times n}) \geq f(m,n)$, where $f(m,n)$ is the desired RHS in the claim. We will now show that $\beta(Q_{m \times n}) \leq f(m,n)$ in all cases.

Let $m\geq n$ WLG. In the trivial case where $m \geq 3n-2$, it's clear that $n$ rows by themselves suffice to cover the $m\times n$ board. Hereafter we consider the non-trivial case.

We refer to the case where $m+n \equiv 2 \mod 4$ as the critical case. Among these, the hard critical case occurs when $m,n$ are even and $m+n \equiv 6 \mod 8$ or when $m,n$ are odd and $m+n \equiv 2 \mod 8$. Boards which are critical but not hard critical constitute the easy critical case. All non-critical boards can be reduced to some easy critical board using Observation \ref{Obs:Monotonicity} by increasing $m$ or $n$ or both without increasing $f(m,n)$. In fact even the hard critical case can be reduced to the easy critical case by increasing both $m$ and $n$ by $2$. It suffices to show that $\beta(Q_{m \times n}) \leq f(m,n)$ in the easy critical case.

Let us now consider the easy critical case. We know that $(m+n-2)/4 < n-1$ since this is equivalent to $m<3n-2$. Let $p=f(m,n)=(m+n-2)/4$ and choose the rightmost $p$ columns and the topmost $p$ rows. The set of uncovered cells now forms a $m' \times n'$ board $V(B_{m'\times n'})$ where $m'=m-p$ and $n'=n-p$ are both even. Lemma \ref{Lemma:BishopsRectangularBoardsFullSolution} implies that  $p'$ sum diagonals and $p'$ difference diagonals can cover $V(B_{m'\times n'})$, where
\begin{align*}
p' &= (m'+n'-2)/2 \\
&=(m-p+n-p-2)/2\\
&= (4p - 2p)/2 = p
\end{align*}
which implies that $\beta(Q_{m \times n}) \leq p = f(m,n)$.
\end{proof}

We've fully solved the Relaxed Queen's Domination problem on rectangular boards.

\section{Queen's Domination of $(4k+3)\times(4k+3)$ boards}

In this section we will show some results and make some conjectures towards the goal of simplifying the long, complex proof by Weakley et al \cite{Weakley1995, Weakley2002, Finozhenok2007} that $\gamma(Q_{n\times n})\geq (n+1)/2$ when $n \equiv 3 \mod 4$ and $n > 11$. This is the only case for square boards in which we haven't simplified the proof of the best-known lower bound, as summarized in Table \ref{Table:SquareBoardsSummary}. For this, we will continue choosing rows, columns and diagonals rather than placing queens. Thus, having solved the Relaxed Queen's Domination problem completely for square (and rectangular) boards, we now investigate the properties of its optimal solutions.

Theorem \ref{Thm:AllSquareBoards} implies that $\gamma(Q_{n\times n})\geq (n-1)/2$ and what remains is to improve the lower bound by $1$ by showing that equality cannot occur. The aforementioned proof spanning 3 papers and over 35 pages did this by proving a sequence of results, stated below. In stating all these results, we assume that $n=4k+3$ is a positive integer such that $\gamma(Q_{n\times n})= (n-1)/2=2k+1$ and $\mathcal D$ denotes the corresponding dominating set of queens. We also consider the sub-board $U$ bounded by the leftmost and rightmost uncovered columns and the bottommost and topmost uncovered rows.

\begin{result}\label{Res:4kPlus3PrelimResult}
   $U,\mathcal D$ have the following properties:
  \begin{enumerate}[label=(\alph*)]
    \item $\mathcal D$ is independent.
    \item $U$ is a square sub-board of size $j \times j$ where $j$ is odd and $3(n+1)/4 \leq j \leq n$.
    \item Each edge square of $U$ is attacked exactly once.
    \item For each queen in $\mathcal D$, either the row or the column through it intersects $U$. Moreover both diagonals through it intersect the edge cells of $U$ in two cells each.
    \item The sum diagonal joining the top-left and bottom-right corner of $U$ covers exactly one queen in $\mathcal D$. Of the remaining $2k$ queens, half are above it and half below it. A similar result holds for the difference diagonal joining the two other corners of $U$.
  \end{enumerate}
\end{result}

We use Result \ref{Res:4kPlus3PrelimResult} to define the center of $U$ as the new origin for our co-ordinate system, similar to \cite{Weakley2002}. Results \ref{Res:4kPlus3HeavyLifting}, \ref{Res:4kPlus3DiagonalNumbers} and Theorem \ref{Thm:UniformSpacedGridPerfectCovers} use this new origin, in contrast to previous sections which used the bottom-left square as the origin. The columns and row numbers still increase from left-to-right and bottom-to-top respectively. Hence, the square $U$ is bounded by the rows and columns numbered $-(j-1)/2$ and $(j-1)/2$.

\begin{result}\label{Res:4kPlus3HeavyLifting}
$U=V(Q_{n\times n})$. Each queen in $\mathcal D$ has both co-ordinates even.
\end{result}

\begin{result}\label{Res:4kPlus3DiagonalNumbers}
There exists some integer $e$, $0\leq e \leq k$ such that the $2k+1$ occupied sum diagonal numbers and difference diagonal numbers are both \[0,\pm 2,\pm 4,\ldots\pm 2e,\pm (2e+4),\pm(2e+8),\ldots \pm (4k-2e)\].
\end{result}

\begin{result}\label{Res:4kPlus3FinalResult}
$n \in \{3,11\}$
\end{result}

Result \ref{Res:4kPlus3PrelimResult} was proven in \cite{Weakley1995} and is equivalently stated as Theorem 2 in \cite{Weakley2002}. As we'll see, the ideas from the proofs of Theorem \ref{Thm:AllSquareBoards} and Theorem \ref{Thm:4k1SquareBoards} suffice to generalize and simplify the proof of most of it. Result \ref{Res:4kPlus3HeavyLifting} is the hardest part of the proof. It spans 19 pages \cite{Weakley2002} and uses 8 lemmas. We propose open questions and make conjectures towards simplifying it. Result \ref{Res:4kPlus3DiagonalNumbers} is the easiest of the four results. We generalize it and provide an alternative proof that is easier to visualize. The proof of result \ref{Res:4kPlus3FinalResult} is based on the Pell's equation \cite{Finozhenok2007, Weakley2002}. This can't be simplified by just choosing occupied rows, columns and diagonals since every solution from Result \ref{Res:4kPlus3DiagonalNumbers} is valid for the Relaxed Queen's Domination problem for all $n \equiv 3 \mod 4$.

For the rest of this section, let $n=4k+3$ be a positive integer. Let $p=\beta(Q_{n \times n})=(n-1)/2=2k+1$ and let rows $(R_i)_{i=0}^{p-1}$, columns $ (C_i)_{i=0}^{p-1}$, sum diagonals $(S_i)_{i=0}^{p-1}$ and difference diagonals $(D_i)_{i=0}^{p-1}$ - collectively referred to as the $4p$ chosen lines - constitute a solution to the Relaxed Queens Domination Problem for a $n\times n$ board. Parts (a) through (e) of the following Theorem \ref{Thm:4kPlus3PrelimResultGeneralized} generalize parts (a) through (e) of Result \ref{Res:4kPlus3PrelimResult} respectively.

\begin{theorem}\label{Thm:4kPlus3PrelimResultGeneralized}
    Consider the sub-board $U$ bounded by the leftmost $(l)$ and rightmost $(r)$ unchosen columns and the bottommost $(b)$ and topmost $(t)$ unchosen rows. Then
  \begin{enumerate}[label=(\alph*)]
    \item The rows $(R_i)_{i=0}^{p-1}$ are all distinct. Similarly, the columns $ (C_i)_{i=0}^{p-1}$, sum diagonals $(S_i)_{i=0}^{p-1}$ and difference diagonals $(D_i)_{i=0}^{p-1}$ are each also all distinct.
    \item $U$ is a square sub-board.
    \item Each edge square of $U$ is covered by exactly one of the $4p$ chosen lines.
    \item Each chosen diagonal intersects the edge cells of $U$ in two points each.
    \item The sum diagonal joining the top-left and bottom-right corner of $U$ is chosen. Of the remaining $2k$ chosen sum diagonals, half are above it and half are below it. An analogous result holds for the difference diagonals.
  \end{enumerate}
\end{theorem}
\begin{proof}
We will re-use the ideas in the proof of Theorem \ref{Thm:AllSquareBoards} and Theorem \ref{Thm:4k1SquareBoards}. Observe that $p=\beta(Q_{n \times n})=2k+1$ implies that equality holds in Equation \ref{eqn:AllSquareBoardsEquation} of Theorem \ref{Thm:AllSquareBoards}.
\begin{enumerate}[label=(\alph*)]
\item For equality to hold in Equation \ref{eqn:AllSquareBoardsEquation}: (i) the rows and columns must both be all distinct for $|K|$ to equal $4(n-p)-4=8k+4$ (ii) the $2k+1$ sum and $2k+1$ difference diagonals must both be all distinct for the $8k+4$ cells in $|K|$ to be covered.
\item The equality in Equation \ref{eqn:AllSquareBoardsEquation} cannot hold if $r-l \neq t-b$, as was shown using the colouring argument in Case 2 of Theorem \ref{Thm:4k1SquareBoards}. Note that the argument, originally stated for integers of the form $4k+1$ also applies to those of the form $4k+3$.
\item Let $U_e$ denote the edge cells of $U$, so that $K \subseteq U_e$. The $8k+4$ cells in $K$ are all covered exactly once, two each by the $4k+2$ chosen diagonals. WLG consider a cell of $U_e \setminus K$ in the top row $t$. It is covered by some chosen column between $l$ and $r$. Moreover, it is not covered by any other chosen line since the diagonal intersections with $U_e$ are all accounted for by $K$ and since $t$ is unchosen by definition.
\item Already established in the proof of part (c).
\item This result was proven in Case 1 of Theorem \ref{Thm:4k1SquareBoards}. Note that the argument, originally stated for integers of the form $4k+1$ also applies to those of the form $4k+3$. It doesn't lead to a contradiction in this case. In fact $\beta(Q_{n\times n})$ equals $(n-1)/2$.
\end{enumerate}
\end{proof}

Most of the subresults (a) through (e) of Result \ref{Res:4kPlus3PrelimResult} are immediately implied by subresults (a) through (e) of Theorem \ref{Thm:4kPlus3PrelimResultGeneralized} respectively.  The latter results generalize the former results to the Relaxed Queen's Domination setting. Theorem \ref{Res:4kPlus3PrelimResult} (b) is the only result that is only partially generalized. Specifically, the parts that are not generalized are two properties of $j$, the length of square $U$: (i) $j$ is odd (ii) $3(n+1)/4 \leq j \leq n$. Note that both these properties aren't true at all in the Relaxed Queen's Domination setting. In fact one can find a counterexample for both these for every $n \equiv 3 \mod 4$, analogous to the counterexample for $n=11$ that is depicted in Figure \ref{figure:SquareBoardsConstruction}.

Having generalized Result \ref{Res:4kPlus3PrelimResult}, we now propose open questions and conjectures towards simplifying Result \ref{Res:4kPlus3HeavyLifting}. Among these, Question \ref{Q:ThirdOpenQ} is particularly simple to state independent of the Queen's Domination Problem and may be of independent interest. We ask what values the set of rows $(R_i)_{i=0}^{p-1}$ and the set of columns $ (C_i)_{i=0}^{p-1}$ can take. This is our first open question and some follow-ups to this are the subsequent open questions.

To state the first open question, we will focus not on the set of chosen rows and columns but the unchosen ones. Let $(R'_i)_{i=0}^p$ denote the unchosen rows and $(C'_i)_{i=0}^p$ denote the unchosen columns. WLG let both these be sorted, so that $0\leq R'_0<\ldots< R'_p<n=4k+3$ and $0\leq C'_0<\ldots< C'_p<n=4k+3$. The set of unchosen cells is equal to the spaced grid $G(\{C'_i\}_{i=0}^{p},\{R'_i\}_{i=0}^{p})$ - shorthand $G$ - as defined in Definition \ref{Def:SpacedGrid}. Theorem \ref{Thm:4kPlus1RectangularBoardsGeneralization} (b) implies that $C'_p-C'_0$ equals $R'_p-R'_0$. WLG, we translate the rows and columns so that $0=R'_0<\ldots<R'_p=e$ and $0=C'_0<\ldots< C'_p=e$. Recall that $p=2k+1$ and that $p$ sum diagonals $(S_i)_{i=0}^{p-1}$ and 
$p$ difference diagonals $(D_i)_{i=0}^{p-1}$ fully cover the $(p+1) \times (p+1)$ spaced grid $G$. Lemma \ref{Lemma:SpacedOddRectangularGrid} implies that $p$ is in fact the least such number that can suffice. Let us now state the first open question, generalizing the $e<4k+3$ constraint: \textit{Which are the `perfectly-coverable' spaced grids with equal number of rows and columns?}

\begin{question}\label{Q:FirstOpenQ}
Let $p=2k+1$ and $e\geq p$. What are the values $0=R'_0<\ldots< R'_p=e$ and $0=C'_0<\ldots< C'_p=e$ such that $G(\{C'_i\}_{i=0}^{p},\{R'_i\}_{i=0}^{p})$ can be fully covered using $p$ sum diagonals and $p$ difference diagonals?
\end{question}

Each of the $2p$ chosen diagonals must cover $2$ points in $K$, the set of outermost points in $G$ which touch the boundary of its convex hull. This is because equality in Lemma \ref{Lemma:SpacedOddRectangularGrid} implies equality in Equation \ref{eqn:SpacedOddRectangularGridEquation}. So far we've only considered the cells in $K\subseteq G$ but we will now consider all cells in $G$. We denote $\{R'_0,\ldots R'_p\}$ by $R'$ and $\{C'_0,\ldots C'_p\}$ by $C'$. Each cell $(C'_i, R'_j)$ in $G$ must be covered by either a sum or a difference diagonal which also covers two cells in $K$. This means that either (i) $C'_i+R'_j \mod e \in R' \cap C'$ or (ii) $C'_i-R'_j \mod e \in C'$ and $R'_j-C'_i \mod e \in R'$. If neither of these is true, then the cell $(C'_i, R'_j)$ is obviously uncoverable. This leads us to our second open question: \textit{Which are the spaced grids with equal number of rows and columns that have no obviously uncoverable cells?}

\begin{question}\label{Q:SecondOpenQ}
Let $e\geq p\geq 2$. What are the values $0=C'_0<\ldots< C'_p=e$ and $0=R'_0<\ldots< R'_p=e$ such that for every $(i,j); 0\leq i,j\leq p$ one of the following holds:
\begin{enumerate}[label=(\roman*)]
\item $C'_i+R'_j \mod e \in C' \cap R'$
\item $C'_i-R'_j \mod e \in C'$ and $R'_j-C'_i \mod e \in R'$
\end{enumerate}
\end{question}

For the Queen's Domination problem for $(4k+3)\times(4k+3)$ boards, we are only interested in odd $p$. We wrote a script\footnote{\url{https://github.com/architkarandikar/queens-domination/tree/main/qstn_5p7_cjtr_5p8}} to investigate this question for small $p$ and $e$ and found that all such spaced grids are symmetric about both their diagonals for odd $p$. At least one of these two symmetries holds for even $p$. This leads us to Conjecture \ref{C:SecondOpenQConjecture} which states some necessary but not sufficient properties for Question \ref{Q:SecondOpenQ}. The first part of this conjecture states the following: \textit{Spaced grids with an equal and odd number of rows and columns that have no obviously uncoverable cells are symmetric about both diagonals.}

\begin{conjecture}\label{C:SecondOpenQConjecture}
All solutions to Question \ref{Q:SecondOpenQ} for odd $p$ satisfy:
\begin{enumerate}[label=(\roman*)]
\item $C'_i=R'_i$ for all $0\leq i \leq p$
\item $C'_i+R'_{p-i}=e$ for all $0\leq i \leq p$
\end{enumerate}
For even $p$, at least one of (i) or (ii) is satisfied.
\end{conjecture}

Conjecture \ref{C:SecondOpenQConjecture} says that solutions to Question \ref{Q:SecondOpenQ} for odd $p$ are symmetric about both diagonals. This leads us to ask a simpler question: does one symmetry imply the other? We already know that the answer is no for even $p$. WLG, we assume symmetry along the sum diagonal, so that $C'_i=R'_i$. The two conditions of Question \ref{Q:SecondOpenQ} now simplify further: for all $C'_i, C'_j\in C'$, one of the following should be true (i) $C'_i+C'_j \mod e\in C'$ (ii) $C'_i-C'_j \mod e\in C'$. The question now becomes: which are such subsets $C'$? Since we've translated so that $C'_0=0$, we are only interested in those $C'$ for which $0\in C'$. For $(4k+3)\times(4k+3)$ boards we are interested in odd-sized subsets. However, we generalize both these constraints to make the question even more natural. The resultant Question \ref{Q:ThirdOpenQ} is now so simple to state that it is potentially of independent interest as a problem. Although it's stated independently of the Queen's Domination problem, it essentially asks: \textit{Which are the spaced grids with equal number of rows and columns symmetric about one diagonal that have no obviously uncoverable cells?}

\begin{question}\label{Q:ThirdOpenQ}
What are the subsets $S$ of $\{0,1,\ldots,e-1\}$ such that for all $x,y\in S$, either $x+y\mod e$ or $x-y\mod e$ is in $S$?
\end{question}

We also wrote a script\footnote{\url{https://github.com/architkarandikar/queens-domination/tree/main/qstn_5p9_cjtr_5p10}} to investigate this question for small values for of $e$. For odd-sized subsets we found the symmetry properties we expected from Conjecture \ref{C:SecondOpenQConjecture}. For even-sized subsets these do not hold. This leads us to Conjecture \ref{C:ThirdOpenQConjecture}, which states necessary but not sufficient properties for Question \ref{Q:ThirdOpenQ}: \textit{Odd-sized subsets which are solutions to Question \ref{Q:ThirdOpenQ} always contain $0$ and are symmetric about $e/2$ except for $0$.}

\begin{conjecture}\label{C:ThirdOpenQConjecture}
All odd-sized solutions $S$ to Question \ref{Q:ThirdOpenQ} satisfy
\begin{enumerate}[label=(\roman*)]
\item $0 \in S$
\item $S \cup \{e\}$ is symmetric about $e/2$
\end{enumerate}
\end{conjecture}

\begin{figure}[t]
\centering
    \begin{subfigure}[t]{.5\linewidth}
        \centering
        \begin{tabular}{c}
\begin{lstlisting}






C':0 1 2 6 7 8  R':0 1 2 6 7 8 
C':0 2 3 5 6 8  R':0 2 3 5 6 8
\end{lstlisting}
        \end{tabular}
        \caption{Solutions to Qn. \ref{Q:SecondOpenQ} for $p=5,e=8$.}
    \end{subfigure}%
    \begin{subfigure}[t]{.5\linewidth}
        \centering
        \begin{tabular}{c}
\begin{lstlisting}
C':0 1 2 5 6  R':0 1 2 5 6 
C':0 1 2 5 6  R':0 1 4 5 6 
C':0 1 3 4 6  R':0 1 3 4 6 
C':0 1 3 4 6  R':0 2 3 5 6 
C':0 1 4 5 6  R':0 1 2 5 6 
C':0 1 4 5 6  R':0 1 4 5 6 
C':0 2 3 5 6  R':0 1 3 4 6 
C':0 2 3 5 6  R':0 2 3 5 6 
\end{lstlisting}
        \end{tabular}
        \caption{Solutions to Qn. \ref{Q:SecondOpenQ} for $p=4,e=6$.}
    \end{subfigure}
    \begin{subfigure}[t]{.5\linewidth}
        \centering
        \begin{tabular}{c}
\begin{lstlisting}









S u {e}: 0 6
S u {e}: 0 2 4 6
S u {e}: 0 1 5 6
S u {e}: 0 1 2 4 5 6
\end{lstlisting}
        \end{tabular}
        \caption{Odd-sized solutions S to Qn. \ref{Q:ThirdOpenQ} for $e=6$.}
    \end{subfigure}%
    \begin{subfigure}[t]{.5\linewidth}
        \centering
        \begin{tabular}{c}
\begin{lstlisting}
S u {e}: 6
S u {e}: 0 1 6
S u {e}: 0 2 6
S u {e}: 0 3 6
S u {e}: 0 4 6
S u {e}: 2 4 6
S u {e}: 0 1 3 4 6
S u {e}: 0 5 6
S u {e}: 0 1 2 5 6
S u {e}: 0 2 3 5 6
S u {e}: 0 1 4 5 6
S u {e}: 1 2 4 5 6
S u {e}: 0 1 2 3 4 5 6
\end{lstlisting}
        \end{tabular}
        \caption{Even-sized solutions S to Qn. \ref{Q:ThirdOpenQ} for $e=6$.}
    \end{subfigure}
\caption{Explanatory examples for Questions \ref{Q:SecondOpenQ}, \ref{Q:ThirdOpenQ} and Conjectures \ref{C:SecondOpenQConjecture}, \ref{C:ThirdOpenQConjecture}.}
\end{figure}

In the last part of this section  we generalize Result \ref{Res:4kPlus3DiagonalNumbers} by characterizing all perfect diagonal covers for Uniformly Spaced Grids with an equal number of rows and columns. 

\begin{definition}[Perfect Diagonal Cover of a Spaced Grid]
A perfect diagonal cover for a spaced grid with $m'$ rows and $n'$ columns is a set of exactly $(m'+n'-2)/2$ sum diagonals and $(m'+n'-2)/2$ difference diagonals that covers every cell in the it.
\end{definition}

Lemma \ref{Lemma:SpacedOddRectangularGrid} implies that a perfect diagonal cover can exist only when $m',n'$ are both even. Theorem \ref{Thm:UniformSpacedGridPerfectCovers} characterizes all perfect diagonal covers for Uniformly Spaced Grids with $p$ rows and $p$ columns when $p=2k+2$ is even. We parameterize $p$ as $2k+2$ because the spaced grid of uncovered cells that remains for $n=4k+3$ after choosing $(n-1)/2$ rows and columns has $2k+2$ unchosen rows and columns. With $m'=n'=p=2k+2$, perfect diagonal covers have $(m'+n'-2)/2=2k+1$ diagonals of each kind.

\begin{theorem}\label{Thm:UniformSpacedGridPerfectCovers}
Consider a Uniformly Spaced Grid $G(\{C'_i\}_{i=0}^{p-1},\{R'_i\}_{i=0}^{p-1})$ with $p=2k+2$ columns $\{C'_i\}_{i=0}^{p-1}$, $p$ rows $\{R'_i\}_{i=0}^{p-1}$ and spacing $d$ between consecutive rows and columns. We identify columns by their $x$-coordinate and rows by their $y$-coordinate. WLG translate the columns and rows to be symmetric about both axes, so that $C'_i=(i-k-1/2)d$ and $R'_i=(i-k-1/2)d$. A set of $2k+1$ sum diagonals $S=\{S_i\}_{i=0}^{2k}$ and $2k+1$ difference diagonals $D=\{D_i\}_{i=0}^{2k}$ - identified by the invariant sum and difference of their cell co-ordinates respectively - is a perfect cover if and only if \begin{gather*}S=D=Q_e \text{ for some } e \text{ such that } 0\leq e\leq k \text{ where }\\
Q_e=\{0,\underbrace{\pm d, \pm 2d, \ldots, \pm ed}_{\text{count: $2e$}}, \underbrace{\pm (e+2)d, \pm (e+4)d, \ldots, \pm (2k-e)d}_{\text{count: $2(k-e)$}}\}\end{gather*}
\end{theorem}
\begin{proof}
WLG we prove this for $d=2$. Note that all $R_i$ and $C_i$ values are odd integers and all $S_i$ and $D_i$ values are even integers for $d=2$. Our proof is based on the following \textit{central observation: all $G$-cells in difference diagonal $2f$ are covered if and only if either difference diagonal $2f$ is selected or if sum diagonals $H_f=\{2f-4k-2,2f-4k+2,\ldots,4k+2-2f\}$ are all selected.} This is because the $G$-cells in difference diagonal $2f$ are $\{(2f-(2k+1),-(2k+1)),(2f-(2k-1),-(2k-1)),\ldots,(2k+1,(2k+1)-2f)\}$. Similar results also hold for difference and sum diagonals $\pm 2f$ for all $f \in [1,2k+1]$.

Let $(S,D)$ be a perfect cover of $G$ and consider the largest integer $e \geq 0$ such that $\{0,\pm 2,\pm 4 \ldots, \pm 2e\} \subseteq S \cap D$. If $e=k$, the claim trivially holds since $|S|=|D|=2k+1$ diagonals of each kind are all already accounted for by $S \cap D$. Otherwise one of $\pm2(e+1)$ does not belong to either $S$ or $D$. WLG, let $2(e+1) \notin D$. The central observation now implies that $H_{e+1}=\{2e-4k,2e-4k+4,\ldots,4k-2e\} \subseteq S$. Hence $\{\pm2(e+2), \pm2(e+4) \ldots, \pm2(2k-e)\} \subseteq S$ implying that $S=Q_e$, since all $2k+1$ elements of $S$ are accounted for. Hence $2(e+1)\notin S$, which similarly implies that $D=Q_e$.

Conversely let $S=D=Q_e$ for some integer $e$ such that $0\leq e\leq k$. We will show that $S \cup D$ covers all cells in $G$ which trivially implies that $(S,D)$ is a perfect cover of $G$. Each cell $a$ in $G$ belongs to some difference diagonal $2f$ where $-2k+1\leq f \leq 2k+1$. Difference diagonals in $Q_e$ are trivially covered since they are selected. Those not in $Q_e$ are $A_e \cup B_e$ where $A_e=\{\pm2(e+1), \pm2(e+3), \ldots, \pm2(2k-e-1)\}$ and $B_e=\{\pm2(2k-e+1), \pm2(2k-e+2), \ldots, \pm2(2k+1)\}$. We infer from the central observation that all cells in difference diagonal $\pm 2(e+1)$ are covered since $H_{e+1}=\{2e-4k,2e-4k+4,\ldots,4k-2e\} \subseteq S$. The central observation further implies that all cells in each difference diagonal in $A_e$ are covered since $H_{2k-e-1}\subset \ldots \subset H_{e+3} \subset H_{e+1}$ and that all cells in each difference diagonal in $B_e$ are covered since $H_i \subseteq \{0,\pm 2,\pm 4 \ldots, \pm 2e\} \subseteq S$ for all integers $i$ such that $2k-e+1\leq i\leq 2k+1$. Hence all cells in $G$ are covered by $S \cup D$.
\end{proof}

\begin{figure}[t]
\includegraphics[width=12cm]{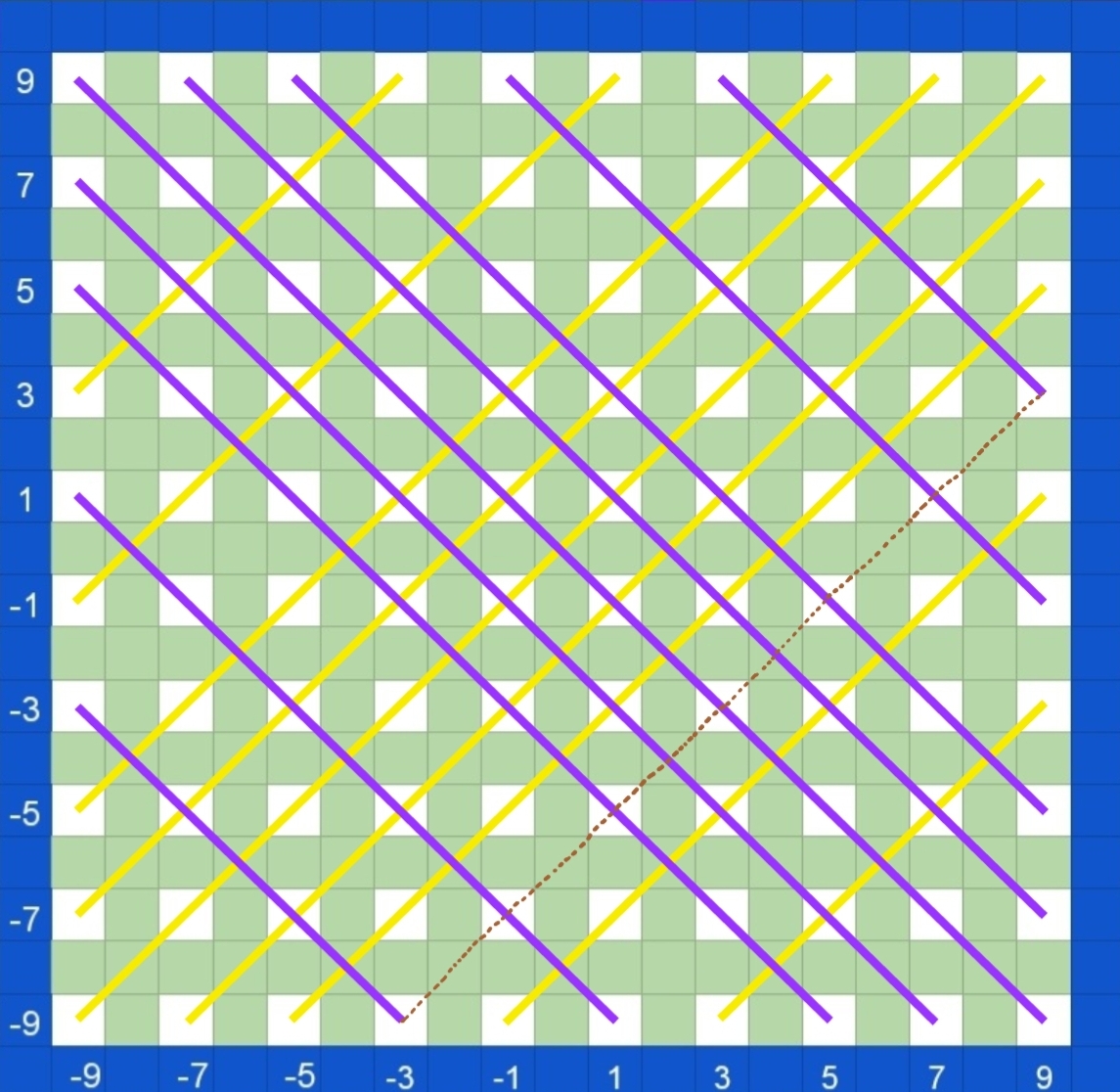}
\centering
\caption{A perfect cover for a Uniformly Spaced Grid with $10$ rows and $10$ columns, which illustrates the proof of Theorem \ref{Thm:UniformSpacedGridPerfectCovers} for $p=10$, $k=4$, $d=2$, $e=2$ and $S=D=Q_e=\{0,\pm 2, \pm 4, \pm 8, \pm 12\}$. Chosen sum and difference diagonals are shown in purple and yellow respectively. The dotted brown difference diagonal $2f=6$ for $f=3$ illustrates the central observation. Since it is unchosen, the sum diagonals $\{-12, -8, -4, 0, 4, 8, 12\}$ must all be chosen so that all cells in it are covered.}
\label{figure:UniformSpacedGrid}
\end{figure}

This proof of Theorem \ref{Thm:UniformSpacedGridPerfectCovers} illustrated in Figure \ref{figure:UniformSpacedGrid} characterizes all perfect covers of a Uniformly Spaced Grid. We have now generalized or made conjectures towards simplifying Result \ref{Res:4kPlus3PrelimResult}, Result \ref{Res:4kPlus3HeavyLifting} and Result \ref{Res:4kPlus3DiagonalNumbers} while only choosing occupied rows, columns and diagonals rather than placing queens. We've noted that Result \ref{Res:4kPlus3FinalResult} can't be simplified in this setting since solutions which choose all even columns, all even rows and diagonals as per Theorem \ref{Thm:UniformSpacedGridPerfectCovers} are all valid for the Relaxed Queen's domination problem. This concludes the last section of this work.

\subsection*{Conclusion}

We've improved the lower bound for the Queen's Domination Problem on rectangular boards by proposing and completely solving the Relaxed Queen's Domination problem. By doing so, we've answered an open question posed in \cite{Bozoki2019} in the affirmative. In the process, we've simplified the proofs for lower bounds on square boards and generalized Weakley's improved lower bound \cite{Weakley1995} for $(4k+1)$-square boards to rectangular boards. Finally we've formulated open questions and conjectures relevant to this problem.

\subsection*{Acknowledgements}

Thanks to our good friend Varun Dwarkanathan for introducing this problem to us.


\nocite{*}

\bibliographystyle{plain} 
\bibliography{refs} 

\begin{thebibliography}{10}

\bibitem{Ahrens1901}
W.~Ahrens.
\newblock {\em Mathematische Unterhaltungen und Spiele}.
\newblock B.G. Teubner, 1901.

\bibitem{Ahrens1910}
W.~Ahrens.
\newblock {\em Mathematische Unterhaltungen und Spiele}.
\newblock B.G. Teubner, 1910.

\bibitem{ball1892mathematical}
W.W.R. Ball.
\newblock {\em Mathematical Recreations and Essays}.
\newblock Read Books, 1892.

\bibitem{Bozoki2019}
Sándor Bozóki, Péter Gál, István Marosi, and William~D. Weakley.
\newblock Domination of the rectangular queen's graph.
\newblock 2019.

\bibitem{BurgerThesis}
A.P Burger.
\newblock {\em The Queens Domination Problem}.
\newblock PhD thesis, University of South Africa, 1998.

\bibitem{Cockayne199013}
E.J. Cockayne.
\newblock Chessboard domination problems.
\newblock {\em Discrete Mathematics}, 86(1):13--20, 1990.

\bibitem{Cockayne1986}
E.J Cockayne, B~Gamble, and B~Shepherd.
\newblock Domination parameters for the bishops graph.
\newblock {\em Discrete Mathematics}, 58(3):221--227, 1986.

\bibitem{Finozhenok2007}
Dmitry Finozhenok and William Weakley.
\newblock An improved lower bound for domination numbers of the queen's graph.
\newblock {\em Mathematical Sciences Faculty Publications}, 37, 01 2007.

\bibitem{UnsolvedProblemNTBook}
Richard~K. Guy.
\newblock {\em Unsolved problems in number theory}.
\newblock Problem Books in Mathematics. Springer-Verlag, New York, second
  edition, 1994.
\newblock Unsolved Problems in Intuitive Mathematics, I.

\bibitem{jaenisch1862}
C.F. Jaenisch.
\newblock {\em Trait{\'e} des applications de l'analyse math{\'e}matique au jeu
  des {\'e}checs, pr{\'e}c{\'e}d{\'e} d'une introduction {\`a} l'usage des
  lecteurs soit {\'e}trangers aux {\'e}checs, soit peu vers{\'e}s dans
  l'analyse}.
\newblock Number v. 1 in Cornell University Library historical math monographs.
  Dufour et cie, 1862.

\bibitem{Knuth2000}
Donald~E. Knuth.
\newblock Dancing links.
\newblock 2000.

\bibitem{stergrd2001ValuesOD}
Patric R.~J. {\"O}sterg{\aa}rd and William~D. Weakley.
\newblock Values of domination numbers of the queen's graph.
\newblock {\em Electron. J. Comb.}, 8, 2001.

\bibitem{Raghavan1987281}
Vijay Raghavan and Shankar~M. Venkatesan.
\newblock On bounds for a board covering problem.
\newblock {\em Information Processing Letters}, 25(5):281--284, 1987.

\bibitem{Szily1902}
K.~von Szily.
\newblock {\em Das Minimalproblem der Damen}, pages 326--328.
\newblock Deutsche Schachzeitung, 57, 1902.

\bibitem{Szily1903}
K.~von Szily.
\newblock {\em Das Minimalproblem der Damen}, pages 65--68.
\newblock Deutsche Schachzeitung, 58, 1903.

\bibitem{Watkins2004}
John~J. Watkins.
\newblock {\em Across the Board: The Mathematics of Chessboard Problems}.
\newblock Princeton University Press, 2004.

\bibitem{Weakley1995}
W.D. Weakley.
\newblock Domination in the queen's graph.
\newblock In Y.~Alavi and A.J.n Schwenk, editors, {\em Graph Theory,
  Combinatorics, and Algorithms}, volume~2, pages 1223--1232. Wiley, Oxford,
  1995.

\bibitem{Weakley2002}
W.D. Weakley.
\newblock A lower bound for domination numbers of the queen's graph.
\newblock In {\em JCMCC: The Journal of Combinatorial Mathematics and
  Combinatorial Computing}, volume~43, pages 231--254. 2002.

\bibitem{Weakley2018}
William~D. Weakley.
\newblock {\em Queens Around the World in Twenty-Five Years}, pages 43--54.
\newblock Springer International Publishing, Cham, 2018.

\bibitem{Weakley2022}
William~D. Weakley.
\newblock Queen domination of even square boards.
\newblock {\em Electron. J. Comb.}, 29(2), 2022.

\end{thebibliography}

\end{document}